\documentclass[12pt,a4paper]{article}
\usepackage[utf8]{inputenc}
\usepackage[T1]{fontenc}
\usepackage{amsmath}
\usepackage{amsfonts}
\usepackage{amsthm}
\usepackage{amssymb}
\usepackage{xcolor}
\usepackage{graphicx}
\usepackage{enumerate}
\usepackage{geometry}
\usepackage{stix}%integrales barrées
\usepackage{hyperref}
\newtheorem{thm}{Theorem}
\newtheorem{rmq}{Remark}
\newtheorem{lemm}{Lemma}

\newtheorem{prop}{Proposition}
\newtheorem{corol}{Corollary}
\newtheorem*{conj}{Conjecture}
\newenvironment{dem}{\noindent\textbf{Proof}}{\qed \\  }

\newcommand{\C}{\mathcal{C}}
\newcommand{\diam}{\mathrm{Diam}}

\newcommand{\s}{\mathbb{S}^2}

\newcommand{\g}{g_{\mathrm{can}}}

\newcommand{\Ric}{\mathrm{Ric}}
\newcommand{\vol}{\mathrm{vol}}

\newcommand{\R}{\mathbb{R}}

\title{Contractive transport maps from $\mathbb{S}^2$ to nearly spherical surfaces with positive Ricci curvature }
\author{Jordan Serres}

\begin{document}
\maketitle

\begin{abstract}
We prove that every nearly spherical, positively curved surface is the contractive, volume-preserving image of a round sphere. The proof combines three main tools: the Ricci flow on surfaces, the Kim–Milman construction, and a multiscale Bakry–Émery criterion.
\end{abstract}

\section{Introduction}

\subsection{Background and main results}

The idea of \emph{comparison geometry} is to relate the properties of a general space, where exact computations are often impossible, to those of certain model spaces whose geometry is well understood. One way for a manifold to be comparable to such model spaces is to have a uniform positive lower bound on the Ricci tensor. Indeed, if $(M^n, g)$ is an $n$-dimensional Riemannian manifold satisfying 
\begin{equation}\label{eq:riccicompar}
\Ric_g \geq (n-1) g,
\end{equation}
then it can be compared to the round sphere $(\mathbb{S}^n, \g)$ of radius $1$, which attains equality in \eqref{eq:riccicompar} and serves as a model space. Many fundamental consequences can be derived, among which:
\begin{itemize}
\item \textbf{(Myers \cite{Meyer})} $\diam(M,g) \leq \diam(\mathbb{S}^n,\g) = \pi$,
\item \textbf{(Bishop-Gromov \cite{bishop1964relation})} $
\vol_g\big(B(x,\rho)\big) \geq r\, \vol_{\g}\big(B(x',\rho)\big),
$ $\forall \rho > 0$, $\forall x \in M^n$, $\forall x' \in \mathbb{S}^n$, 
with the volume ratio $r = \vol_g(M^n) / \vol_{\g}(\mathbb{S}^n)$,
\item \textbf{(Lichnerowicz \cite{Lichnerowicz1958GomtrieDG})} $\lambda_1(M^n, g) \geq \lambda_1(\mathbb{S}^n, \g) = n$, with $\lambda_1$ the first non-zero eigenvalue of the Laplace-Beltrami operator,
\end{itemize}
and also heat-kernel estimates, the Lévy–Gromov isoperimetric comparison, or log-Sobolev inequalities, to name just a few.
Milman \cite{Milman2018} pointed out that all these results can be unified if there exists a map
\[
T: (\mathbb{S}^n, \g, \vol_{\g}) \to (M^n, g, \vol_g)
\]
that is \emph{contractive}, \textit{i.e.}, $1$-Lipschitz, and pushes forward the rescaled volume measure $r\,\vol_{\g}$ onto $\vol_g$, where
\[
r = \frac{\vol_g(M^n)}{\vol_{\g}(\mathbb{S}^n)}.
\]
Indeed, the contraction property immediately implies Myers’s theorem, and for instance, the Lichnerowicz inequality follows from the variational characterization of $\lambda_1$ by noting that
\begin{align*}
n \intbar_{M^n} \left(f - \intbar_{M^n} f\,d\vol_g\right)^2 d\vol_g
&= n \intbar_{\mathbb{S}^n} \left(f \circ T - \intbar_{\mathbb{S}^n} f \circ T\,d\vol_{\g}\right)^2 d\vol_{\g} \\
&\leq \intbar_{\mathbb{S}^n} \left| \nabla \left( f \circ T \right) \right|^2 d\vol_{\g} \\
&\leq \intbar_{\mathbb{S}^n} \left| (\nabla f) \circ T \right|^2 \left|\nabla T\right|^2 d\vol_{\g} \\
&\leq \intbar_{\mathbb{S}^n} \left| (\nabla f) \circ T \right|^2 d\vol_{\g} \\
&= \intbar_{M^n} \left| \nabla f \right|^2 d\vol_g,
\end{align*}
from which it follows that $\lambda_1(M^n, g) \geq n$.
This led Milman to formulate the following conjecture. Note that it is restricted to manifolds diffeomorphic to $\mathbb{S}^n$ so as to avoid obvious topological obstructions to the existence of a Lipschitz continuous map.
\begin{conj}{\cite[Conjecture 4]{Milman2018}}
For any $(\mathbb{S}^n, g, \vol_g)$ satisfying $\Ric_g \geq (n-1) g$, there exists a map from the round metric $\g$ of radius $1$
\[
T : (\mathbb{S}^n, \g, \vol_{\g}) \to (\mathbb{S}^n, g, \vol_g),
\]
which pushes forward $\vol_{\g}$ onto $\vol_g$ up to a finite constant and contracts the corresponding metrics.
\end{conj}
  
The goal of this note is to take a first step toward Milman’s conjecture in dimension $2$ by showing that it holds for nearly spherical surfaces, see Section \ref{sec:mainresults} for a precise statement.

\begin{thm}
If $g$ is a small perturbation of a round metric on $\s$ satisfying $\Ric_g \geq g$, then there exists a contraction \mbox{$T : (\s,\g) \to (\s,g)$} pushing forward the rescaled volume measure $(\vol(\s,g)/4\pi)\,\vol_{\g}$ of the round sphere $\g$ of radius $1$ onto the volume measure $\vol_{g}$.
\end{thm}

The idea of the proof is as follows.  
First, we start a rescaled Ricci flow from the initial metric $g$. This flow converges to the round sphere $\mathbb{S}_\rho^2$ of radius $\rho$, where $\rho^2 = \vol_g(\mathbb{S}) / 4\pi$. Since $\Ric_g \geq g$ implies $\vol_g(\mathbb{S}) \leq 4\pi$, with equality if and only if $g$ is exactly the round sphere of radius $1$, we may assume that $\rho < 1$.  

Second, we apply the Kim-Milman construction \cite{kim2012generalization} to build a transport map from the volume measure of $\mathbb{S}_\rho^2$ onto $\vol_g$. The key point here is to identify the driving vector field of the flow of volume measures, which turns out to be the gradient of the curvature potential, see Lemma \ref{lem:divergencecontinuity}.  

Third, we aim to bound the Lipschitz constant of the Kim-Milman map. This reduces to controlling the Hessian of the curvature potential via a multiscale Bakry-Émery analysis, see Corollary \ref{corol:Bakry-Emerymultiscale}. To this end, we establish a uniform version of Hamilton’s preservation of positive curvature theorem, which can be viewed as a compactness result, see Lemma \ref{lem:Hamiltonuniform}.  

Finally, we obtain a Lipschitz bound for the Kim-Milman map with source $\mathbb{S}_\rho^2$, and the fact that $\rho < 1$ provides just enough room to ensure that composition with the contraction $\s_1 \to \mathbb{S}_\rho^2$ is $1$-Lipschitz whenever the initial metric is sufficiently close to $\mathbb{S}_\rho^2$, see Section \ref{sec:finalproof}.\\

Let us note that dimension $2$ is really used only twice: via the Gauss-Bonnet theorem (see Section \ref{sec:rescaledricci}) and via Hamilton’s theorem (see Lemma \ref{lem:Hamiltonuniform}). The first point does not appear to be essential, as it could be replaced by the monotonicity of the total scalar curvature. The second is more technically involved, but one might try to circumvent it by assuming that the initial manifold $(\mathbb{S}^n, g)$ is conformally equivalent to the round sphere and by using the Yamabe flow. Since the Ricci flow reduces to the Yamabe flow in dimension $2$, this approach would be the most natural first attempt to generalize Theorem \ref{thm:main}.

\subsection{Related literature}

The most famous related result is Caffarelli's contraction theorem \cite{caffarelli1}, stating that if $\mu = e^{-V(x)}dx$ is a probability distribution on $\R^n$ such that $\nabla^2 V\geq I_n$, then the Brenier optimal transport map from the Gaussian $(2\pi)^{-n/2} e^{-|x|^2/2} dx$ onto $\mu$ is a contraction. Its original proof is based on a PDE analysis of the Monge-Ampère equation solved by the Brenier map, but there is now also a proof by means of entropic regularization techniques \cite{caffarellientropy, chewi2023entropic}.

The importance of realizing a probability measure $\mu$ as a Lipschitz image of the Gaussian, not necessarily contractive, has become increasingly recognized in the community \cite{Milman2009}, motivating the search for Lipschitz transport maps beyond the realm of optimal transport. This culminated in the work of Kim and Milman \cite{kim2012generalization}, who introduced a construction of transport maps via heat flows. Their approach provides a precise criterion for when the resulting map, now known as the Kim-Milman map, is Lipschitz continuous, and allows for explicit bounds on its Lipschitz constant. In particular, this method enables the construction of Lipschitz transport maps between broader classes of distributions, notably log-concave measures with certain symmetry assumptions. This seminal construction has been widely applied, see for example \cite{dai2023lipschitz, mikulincer2023lipschitz, brigati2024heat, stephanovitch2024smooth, lopez2025bakry}. Moreover, in \cite{shenfeld2024exact}, it was shown that the Lipschitz continuity of the Kim-Milman map can be derived using the multiscale Bakry-Émery analysis introduced in \cite{bauerbodi} in the context of Euclidean Field Theory. We note the similarity between these ideas and the multiscale Bakry-Émery criterion discussed in Section \ref{sec:multiBE}.
Although originally developed in the Euclidean setting, the Kim-Milman construction can also be adapted to the setting of a heat flow on a Riemannian manifold, as carried out in \cite{fathi2024transportation} for log-Lipschitz perturbations of the volume measure on the round sphere.
In this note, we adapt the Kim-Milman construction to the case of a time-evolving metric on the sphere, see Section \ref{sec:KimMilman}.

In view of Caffarelli’s contraction theorem, Milman’s conjecture appears particularly natural. However, only a few works in the literature have addressed this problem, often focusing on obstructions to related questions \cite{beck2021friedland, fathi2024some}. The present result seems to be the first positive step toward a proof of Milman’s conjecture. The main idea is to use the Ricci flow on positively curved surfaces to construct a transport map via the Kim–Milman approach.

The Ricci flow, introduced by Hamilton in \cite{hamilton1982three} for three-dimensional manifolds, evolves a Riemannian metric according to a nonlinear heat-type PDE that tends to homogenize curvature, and was thus expected to converge toward metrics of constant Ricci curvature. This flow has proved to be remarkably fruitful, breaking down the barriers between topology, differential geometry, PDE theory, and geometry, and culminating in Perelman’s proof of Thurston’s geometrization conjecture \cite{perelman2002entropy, perelman2003ricci}. For detailed expositions, see \cite{chow2006hamilton, kleiner2008notes, ajm/1154098947}.

Hamilton also used the Ricci flow to provide a new proof of the uniformization theorem for closed Riemannian surfaces \cite{hamilton1988ricci}, with the case of positive curvature being the most challenging. Following Hamilton’s work, several simplified alternative proofs were proposed \cite{chow1991ricci, bartz1994new, chen2006note}. However, since our goal is to bound the Hessian of the curvature potential (see Section \ref{sec:unifHamilton}), this note primarily follows Hamilton’s original approach, which is more refined and thus better suited to our purposes.

The literature on the Ricci flow is vast and cannot be summarized in just a few lines. Let us nevertheless mention the differentiable sphere theorem \cite{brendle2010ricci}, the study of Ricci gradient solitons \cite{Naber+2010+125+153}, as well as its place within the broader field of metric flows, which includes, among others, the mean curvature flow \cite{huisken1984flow} and the Yamabe flow \cite{ye1994global}.

\section{Main results and proofs}\label{sec:mainresults}

Our main result can be stated rigorously as follows, using the notation $||\cdot||_{\C^{k,\alpha}}$ for the $k$-th Hölder norm with exponent $\alpha$.

\begin{thm}\label{thm:main}
Let $\alpha\in(0,1)$, and $v\in (0,4\pi)$.
There exists a constant $\varepsilon_0>0$ such that for all smooth Riemannian metric $g$ of volume $v$ on the sphere $\s$ satisfying
\[
|| g - \mathfrak{g} ||_{\C^{5,\alpha}(\s)} \leq \varepsilon_0
\]
with $\mathfrak{g}$ the round metric with the same volume as $g$,
there exists a contraction \mbox{$T: (\s,\g) \to (\s,g)$} pushing forward the rescaled volume measure $(\vol(\s,g)/4\pi)\,\vol_{\g}$ of the round sphere $\g$ of radius $1$ onto the Riemann volume $\vol_{g}$.
\end{thm}

Let us point out that $\varepsilon_0$ depends only on the volume $v=\vol(\s,g)$, and the bound obtained from our proof degenerates to $0$ as $v \to 4\pi$. This behavior is natural, since one cannot expect a general perturbation (even arbitrarily small) of the unit sphere to be a contraction of the unit sphere itself.

The requirement of the $\C^{5,\alpha}$ norm originates from several steps in the argument. Asking the scalar curvature of $g$ to be close to a constant translates into a $\C^{2,\alpha}$ control at the level of the metric. In order to apply the Arzelà–Ascoli theorem in Lemma \ref{lem:Hamiltonuniform}, we further require the gradient of the scalar curvature to be close to zero, which translates into a $\C^{3,\alpha}$ control. Finally, the passage to the $\C^{5,\alpha}$ norm comes from \cite[Theorem A]{bahuaud2020convergence}, again in the proof of Lemma \ref{lem:Hamiltonuniform}, where a $\C^{k+2,\alpha}$ bound on the initial data is needed to obtain a $\C^{k,\alpha}$ bound along the Ricci flow. This stronger assumption thus seems mostly technical, dictated by the method of proof, although sharpening it appears difficult.

Moreover, note that the nearly spherical condition implies $\Ric_g \geq g$, even if it means reducing $\varepsilon_0$ a little more, see Remark \ref{rmq:volumenormalization}. As a consequence, we have the following immediate corollary.

\begin{corol}
Let $g$ be a Riemannian metric on $\s$ such that $\Ric_g\geq g$. Then there exist $\varepsilon_0>0$, only depending on the volume $\vol(\s,g)$, such that if $ || g - \mathfrak{g} ||_{\C^{5,\alpha}(\s)} \leq \varepsilon_0 $,
with $\mathfrak{g}$ the round metric with the same volume as $g$,
then Milman's conjecture holds for $g$.
\end{corol}
 
The remainder of the paper is devoted to the proof of Theorem \ref{thm:main}.

\subsection{The rescaled Ricci flow}\label{sec:rescaledricci}

We start a rescaled Ricci flow $(g_t)_{t\geq 0}$ from the initial metric $g$, in the following way
\begin{equation}\label{eq:rescaledRicciflow}
\left\{
    \begin{array}{ll}
        \partial_t g_t &=  -2 \Ric_{g_t} + r_t g_t \\
        g_0 &=  g
    \end{array}
\right.
\end{equation}
where
\[
r_t :=\intbar_{\s} R_t\,d\nu_t := \frac{\int_{\s} R_t\,d\nu_t }{\int_{\s} d\nu_t}
\]
where $R_t$ stands for the scalar curvature, and $\nu_t$ is the Riemann volume measure for the metric $g_t$. The introduction of the total scalar curvature $r_t$ in the Ricci flow allows to guarantee the conservation of the mass of the Riemann volume along the flow.
And indeed, an immediate calculation, see \textit{e.g.} \cite{topping2006lectures}, shows that $(\nu_t)_t$ evolves according to the following equation
\begin{equation}\label{eq:continuityeq}
\partial_t \nu_t = \left(r_t -R_t\right) d\nu_t 
\end{equation}
In particular, one can check that the mass is conserved, since
\[
\partial_t \int_{\s} d\nu_t = \int_{\s} \partial_t\nu_t = \int_{\s} R_t\,d\nu_t - \int_{\s} R_t\,d\nu_t = 0. 
\]
Therefore the Gauss-Bonnet theorem implies that
\[
r_t = \frac{\int_{\s} R_t\,d\nu_t }{\int_{\s} d\nu_t}  = \frac{8\pi}{\nu_0(\s)} 
\]
since the volume is preserved. Thus, $r_t$ is constant and will be denoted simply by $r$. Note also that we will often omit the time dependence for other quantities, for example, writing $R$ for the scalar curvature, while the context will make it clear that they still depend on time.

If the initial metric $g_0$ has positive scalar curvature, it is well known \cite{hamilton1988ricci, chow1991ricci, bartz1994new} that the rescaled Ricci flow converges exponentially towards the metric $\mathfrak{g}$ of the round sphere of radius $\rho=\sqrt{\nu_0(\s)/4\pi}$, and in particular there is convergence of the volume measures $$\nu_t \underset{t\to\infty}{\longrightarrow} \vol_{\mathfrak{g}}. $$

\subsection{The Kim-Milman map}\label{sec:KimMilman}

To implement the Kim--Milman construction \cite{kim2012generalization} of a transport map from $\vol_{\mathfrak{g}}$ onto $\nu_0 = \vol_g$, we need to express the evolution equation \eqref{eq:continuityeq} in divergence form. This is accomplished by the following lemma.
\begin{lemm}\label{lem:divergencecontinuity}
There exists a familly of functions $\xi_t : \s\to\R$ such that the flow of measures $(\nu_t)_{t\geq 0}$ evolves along the rescaled Ricci flow \eqref{eq:rescaledRicciflow} according to
\begin{equation}\label{eq:divergencecontinuity}
\partial_t \nu_t = -\nabla\cdot \left(\nabla \xi_t\cdot \nu_t \right)
\end{equation}
\begin{rmq}
Let us point out that \eqref{eq:divergencecontinuity} is understood in a weak sense, see \textit{e.g.} \cite[Definition 4.1]{santambrogio2015optimal}, and means that for all smooth enough functions $f:\s\to\R$, it holds that
\[
\partial_t \int_{\s} f\,d\nu_t = \int_{\s} \nabla f\cdot \nabla \xi_t\,d\nu_t.
\]
\end{rmq}
\end{lemm}
\begin{dem}
From \eqref{eq:continuityeq}, we know that 
\[
\partial_t \int_{\s} f\,d\nu_t = \int_{\s} f\left( r -R_t\right) d\nu_t.
\]
So \eqref{eq:divergencecontinuity} holds, if, and only if, for all $f$,
\[
\int_{\s} f\left( r -R_t\right) d\nu_t = \int_{\s} \nabla f\cdot \nabla \xi_t\,d\nu_t.
\]
By using the integration by parts formula, we get that this is equivalent to 
\[
\Delta \xi_t = R_t - r,
\]
which is the Poisson equation in the metric $g_t$ with source term $R_t - r$. Since this source is centered, \textit{i.e}, $\int_{\s} \left(R_t - r\right) d\nu_t = 0$, the Poisson equation admits a unique solution up to a constant \cite[Thm 4.7]{aubin2012nonlinear}. Hence, the existence of the functions $\xi_t$ satisfying \eqref{eq:divergencecontinuity} is ensured.
\end{dem}

\begin{rmq}
The function $\xi_t$, defined as the solution of $\Delta \xi_t = R_t - r$, is called \emph{the curvature potential}. It was used by Hamilton himself in \cite{hamilton1988ricci} in his proof of the convergence of the flow in the case of positive scalar curvature.
\end{rmq}

Now that a divergence form is established for the continuity equation, the Kim-Milman construction can be carried out. We sketch it for clarity purpose, but we refer to \cite{kim2012generalization} for details. It consists in constructing the flow of diffeomorphisms $(S_t)_t$ of the advection field $\nabla \xi_t$ by solving
\[
\partial_t S_t(x) = \nabla \xi_t (S_t(x)),\quad S_0(x)=x.
\]
It is well known that at time $t$, the map $S_t :\s\to\s$ pushes forward the measure $\nu_t$ onto the measure $\nu_0$, see \cite[Theorem 4.4]{santambrogio2015optimal}. The Kim-Milman map is then given by
\begin{equation}\label{eq:defT}
T := \lim_{t\to\infty} S_t^{-1}
\end{equation}
which is a transport from $\nu_\infty=\vol_{\mathfrak{g}}$ onto $\nu_0$.
The goal is then to show that the transport map 
\[
T : (\mathbb{S}, \mathfrak{g},\vol_{\mathfrak{g}}) \to (\mathbb{S}, g_0, \nu_0)
\]
defined by \eqref{eq:defT} is Lipschitz, with a constant sufficiently small to be absorbed by the volume ratio $\nu_0(\mathbb{S}) / 4\pi < 1$. According to the usual Kim-Milman construction, this property can be analyzed via the derivative of the driving vector field $\nabla \xi_t$, \textit{i.e.}, the Hessian $\nabla^2 \xi_t$ of $\xi_t$.

\subsection{A multiscale Bakry-Émery criterion}\label{sec:multiBE}

In this section, we first present a general proposition that ensures the Kim-Milman map, constructed from a flow of Riemannian metrics, is Lipschitz. This result applies to any smooth flow in any dimension. We then derive a corollary for the Ricci flow in dimension $2$, which corresponds to a multiscale Bakry-Émery criterion.

\begin{prop}\label{prop:condipourcontract}
Let $(g_t)_{t\geq 0}$ be a smooth flow of Riemannian metric that converges towards a metric $g_\infty$ on a differential manifold $M$. Assume that the volume is preserved along the flow, and that the volume measures satisfies the weak continuity equation
\[
\partial_t \nu_t = -\nabla\cdot \left(\nabla \xi_t\cdot \nu_t \right)
\]
for some potential functions $\xi_t : M \to \R$. If there exists some $\dot{\lambda}_t \in \R$ such that, for all $t \geq 0$,
\[
2 \nabla^2 \xi_t -\partial_t g_t  \geq -\dot{\lambda}_t\, g_t,
\]
then the Kim-Milman transport map $$ T : (M,g_\infty) \to (M, g_0)$$ defined from the flow is $e^{\int_0^\infty \dot{\lambda}_s\,ds}$-Lipschitz.
\end{prop}

\begin{rmq}
Note that as in Section \ref{sec:KimMilman}, we have that the potential functions $\xi_t$ satisfy the Poisson equation $$\Delta_{g_t} \xi_t = -\frac{1}{2}\mathrm{tr}_{g_t}\,\partial_t g_t.  $$
\end{rmq}

\begin{dem}
First, let us compute the time derivative of $\nabla S_t$. Recall that this gradient is computed with respect to the Riemann metric $g_t$, and therefore also depends on time. Recall also that in charts, it is is given by $$\nabla S_t = g_t^{-1} \nabla_{eucl} S_t, $$ where $\nabla_{eucl} S_t = (\partial_{x_i} S_t^j)_{i,j} $ denotes the euclidean gradient (which is a matrix since $S_t$ is a tangent vector.) We can then compute
\begin{align*}
\partial_t \nabla S_t &= \partial_t(g_t^{-1})\nabla_{eucl} S_t + g_t^{-1}\left(\nabla_{eucl} \partial_t S_t \right) \\
&= - g_t^{-1}(\partial_t g_t) g_t^{-1} \nabla_{eucl} S_t +\nabla \left(\nabla \xi_t (S_t)\right) \\
&= - g_t^{-1}(\partial_t g_t) (\nabla S_t) + \nabla^2\xi_t (\nabla S_t).
\end{align*}
Let us also compute that for all $x\in M$ and all $u\in T_x M$, we have
\begin{align*}
\partial_t |\nabla S_t(x) u |_{g_t}^2 &= |\nabla S_t(x) u |_{\partial g_t}^2 + 2\nabla^2\xi_t \left(\nabla S_t(x) u,\nabla S_t(x) u\right) - 2 \nabla S_t(x) u (\partial_t g_t) \nabla S_t(x) u \\
&= 2\nabla^2\xi_t \left(\nabla S_t(x) u,\nabla S_t(x) u\right) - |\nabla S_t(x) u |_{\partial g_t}^2\\
&= \left( 2\nabla^2\xi_t - \partial_t g_t \right) \left(\nabla S_t(x) u,\nabla S_t(x) u\right).
\end{align*}
Therefore, from the assumption, we get that 
\[
\partial_t |\nabla S_t(x) u |_{g_t}^2 \geq -\dot{\lambda}_t |\nabla S_t(x) u |_{g_t}^2.
\]
So Grönwall's lemma gives
\begin{equation*}\label{eq:calcinterm}
|\nabla S_t(x) u |_{g_t}^2 \geq e^{-\int_0^t \dot{\lambda}_s\,ds} |\nabla S_0(x) u |_{g_0}^2 = e^{-\int_0^t \dot{\lambda}_s\,ds} | u |_{g_0}^2 
\end{equation*}
because by construction of the flow $S_t$, we have that $\nabla S_0 = Id$. Hence the map $$S_t^{-1} : (M,g_t) \to (M, g_0) $$ is $e^{\int_0^t \dot{\lambda}_s\,ds}$-Lipschitz, and letting $t$ go to $+\infty	$ completes the proof.

\end{dem}

Applying Proposition \ref{prop:condipourcontract} in the case of the rescaled Ricci flow \eqref{eq:rescaledRicciflow} on $\s$, we deduce the following multiscale Bakry-Émery criterion, see Section \ref{sec:rescaledricci} for the convergence of the flow.

\begin{corol}\label{corol:Bakry-Emerymultiscale}
Let $(g_t)_{t\geq 0}$ be a rescaled Ricci flow starting from a metric $g_0$ with positive scalar curvature. It converges towards the round metric $\mathfrak{g}$ on $\s$ with same volume as $g_0$. Let $\xi_t : \s \to \R$ be the curvature potential given in Lemma \ref{lem:divergencecontinuity}. If there exists some $\dot{\lambda}_t \in \R$ such that, for all $t \geq 0$,
\[
\nabla^2 \xi_t + \Ric_{g_t}  \geq \frac{1}{2}\left(r -\dot{\lambda}_t\right)\, g_t,
\]
then the Kim-Milman transport map\, $ T : (M,\mathfrak{g}) \to (M, g_0)\,$ is $\,e^{\int_0^\infty \dot{\lambda}_s\,ds}$-Lipschitz.
\end{corol}

\subsection{A uniform version of Hamilton’s theorem}\label{sec:unifHamilton}

In \cite[Theorem 8.1]{hamilton1988ricci}, Hamilton proved that if the initial metric has positive scalar curvature, then there exists a constant $c > 0$, depending only on the initial metric, such that $R \geq c$ along the rescaled Ricci flow on a surface. In this section, we aim to show that the constant $c > 0$ can in fact be chosen uniformly when the initial metric is restricted to the class of nearly round metrics.

\begin{lemm}\label{lem:Hamiltonuniform}
Let $v\in(0,4\pi)$ be a fixed volume, and $\alpha\in (0,1)$ a Hölder exponent. Let $r:=8\pi/v>2$ and $\mathfrak{g}$ denotes the round metric with constant scalar curvature equals to $r$. Let $\varepsilon>0$ small enough such that $r-\varepsilon\geq 2$. Then there exist a constant $C>0$ depending only on $v,\varepsilon$ and $\alpha$, such that all metrics $g$ on $\s$ such that 
\begin{equation}\label{eq:nearlyspherical}
|| g - \mathfrak{g} ||_{\C^{5,\alpha}(\s)} \leq \varepsilon
\end{equation}
satisfy 
\[
\forall t\geq 0,\quad R_t\geq C>0,
\]
where $R_t$ denotes the scalar curvature at time $t$ of the rescaled Ricci flow starting from $g$.
\end{lemm}
\begin{rmq}\label{rmq:volumenormalization}
The assumption $r-\varepsilon\geq 2$ combined with \eqref{eq:nearlyspherical} give that the scalar curvature of $g$ satisfies $R\geq 2$. So this assumption prevents the scalar curvature to be too close to zero.
\end{rmq}
\begin{rmq}\label{rmq:Cepsilongdecroit}
Let us point out that for a fixed $\alpha\in (0,1)$ and $v\in(0,4\pi)$, it is clear that the constant $C=C(\varepsilon)$ of Lemma \ref{lem:Hamiltonuniform} is non-increasing in $\varepsilon$, \textit{i.e.}, $0<\varepsilon\leq\varepsilon' \Rightarrow C(\varepsilon')\leq C(\varepsilon)$.
\end{rmq}

\begin{dem}
The proof proceeds by a contradiction argument. Let us assume that there is a sequence of metrics $g^{(n)}$ satisfying $|| g^{(n)} - \mathfrak{g} ||_{\C^{5,\alpha}} \leq \varepsilon$ and 
\begin{equation}\label{eq:contradictioncompact}
\exists (t_n,x_n)\in\R_+\times\s,\quad R_{t_n}^{(n)}(x_n) \underset{n\to+\infty}{\longrightarrow} 0
\end{equation}
Thanks to Assumption \eqref{eq:nearlyspherical}, the Arzela-Ascoli theorem allows us to extract a subsequence, still denoted $g^{(n)}$, which converges in $\C^{4,\alpha}(\s)$ towards some $g^{(\infty)}$ satisfying $R^{(\infty)}\geq 2$. Now, we use the property of continuous dependance in the initial metric of the Ricci flow, that is for all $\delta>0$ there exists $\delta'>0$ such that if a metric $g$ on $\s$ satisfies $||g^{(\infty)} - g||_{\C^{4,\alpha}(\s)}\leq \delta' $, then it holds that
\begin{equation}\label{eq:conitnuousdependance}
\forall t\geq 0,\quad ||R_t - R_t^{(\infty)}||_\infty \leq  \delta
\end{equation}
A proof of this fact can be found in \cite{bahuaud2020convergence} for the unrescaled Ricci flow, but since the rescaled flow is only a rescaling in space-time of the unrescaled one, the proof can be adapted.
Now there is two cases to distinguish. Either the sequence of times $(t_n)_n$ is bounded or unbouded. If it is bounded, we can extract again subsequences such that $(t_n,x_n)\to (t,x)\in \R_+\times\s$.
Applying \eqref{eq:conitnuousdependance} with $t_n,x_n$ and $R^{(n)}$ and using that $g^{(n)} \overset{\C^{4,\alpha}}{\underset{n\to\infty}{\longrightarrow}} g^{(\infty)} $, we get that $R_{t_n}^{(n)}(x_n) \underset{n\to\infty}{\longrightarrow} R_t^{(\infty)}(x)$. But \eqref{eq:contradictioncompact} gives then that $R_t^{(\infty)}(x)=0$ which contradict the fact that positive scalar curvature is preserved along the rescaled Ricci flow on surfaces, see \cite[Theorem 3.5 and 8.1]{hamilton1988ricci}. Analogously, if $(t_n)_n$ is unbouded, we can extract a subsequence such that $(t_n,x_n)\to (\infty,x)$. Again, applying \eqref{eq:conitnuousdependance} with $t_n, x_n$ and using the convergence of $R^{(n)}$, we obtain  
\[
R_{t_n}^{(n)}(x_n) \underset{n\to\infty}{\longrightarrow} r_\infty,
\]  
where $r_\infty$ is the constant scalar curvature of the round sphere with the same volume as $g^{(\infty)}$, towards which the rescaled Ricci flow starting from $g^{(\infty)}$ converges. And again \eqref{eq:contradictioncompact} gives $r_\infty=0$, which contradicts Hamilton's theorem. The proof is complete.

\end{dem}

\subsection{Proof of Theorem \ref{thm:main}}\label{sec:finalproof}

We can now prove Theorem \ref{thm:main}. Let $\varepsilon>0$ such that \eqref{eq:nearlyspherical} holds.
First, from a direct computation (see \cite[Section 9.3]{hamilton1988ricci}), we get that the curvature potential $\xi_t$, defined in Lemma \ref{lem:divergencecontinuity}, satisfies the following evolution equation along the rescaled Ricci flow \eqref{eq:rescaledRicciflow}
\[
\partial_t |\nabla^2 \xi_t |^2 \leq \Delta (|\nabla^2 \xi_t |^2) -2R|\nabla^2 \xi_t |^2
\]
Using Lemma \ref{lem:Hamiltonuniform}, we have that $R\geq C>0$ for some uniform constant only depending on $\varepsilon$ defined in Assumption \eqref{eq:nearlyspherical}, therefore the maximum principle (see \textit{e.g.} \cite[Theorem 3.1.1]{topping2006lectures}), gives that
\[
\underset{x\in\s}{\sup}|\nabla^2 \xi_t(x) |_{g_t}^2 \leq \underset{x\in\s}{\sup}\,|\nabla^2 \xi_0(x) |_{g_0}^2\, e^{-2Ct}.
\]
In a similar way, setting $M_t := \nabla^2\xi_t - \frac{1}{2}(R-r)g_t$, we have that
\[
\partial |M_t|^2 \leq \Delta (|M_t|^2) -2R|M_t|^2
\]
and therefore again from Lemma \ref{lem:Hamiltonuniform} and the maximum principle, we have that
\[
\underset{x\in\s}{\sup}\, \left|M_t(x)\right|_{g_t}^2 \leq \underset{x\in\s}{\sup}\,|M_0(x)|_{g_0}^2 e^{-2Ct}.
\]
Using the triangle inequality, and denoting $|\nabla^2 \xi_0 |_\infty := \underset{x\in\s}{\sup}\,|\nabla^2 \xi_0(x) |_{g_0}$ and $|M_0|_\infty := \underset{x\in\s}{\sup}\,|M_0(x)|_{g_0} $, we get that
\[
|(R-r)g_t| \leq \left|(R-r)g_t - 2\nabla^2\xi_t  \right| + \left|2\nabla^2\xi_t\right| \leq 4 \max\left\{|\nabla^2 \xi_0 |_\infty\,;\,|M_0|_\infty  \right\} e^{-Ct}
\]
and therefore
\[
\left|2\nabla^2\xi_t + (R-r)g_t \right| \leq 6 \Lambda e^{-C t}
\]
with $\Lambda := \max\left\{|\nabla^2 \xi_0 |_\infty\,;\,|M_0|_\infty  \right\}$.
We deduce that
\[
2\nabla^2\xi_t + (R-r)g_t \geq -6 \Lambda e^{-C t},
\]
from which Corollary \ref{corol:Bakry-Emerymultiscale} gives that the Kim-Milman map from the round sphere of radius $\sqrt{\nu_0(\s)/4\pi}$ onto $(\s,g_0)$ is $e^{3\Lambda/C}$-Lipschitz.  By composing with a dilation, it is therefore a $\sqrt{\nu_0(\s)/4\pi} e^{3\Lambda/C}$-Lipschitz map that pushes forward the rescaled volume measure $(\vol(\s,g)/4\pi)\,\vol_{\g}$ of the round sphere $\g$ of radius $1$ onto the volume measure $\vol_{g_0}$. Consequently, the Kim-Milman map composed with the dilation is a contraction if
\begin{equation}\label{eq:diminuerepsilonpourcontract}
\Lambda \leq \frac{C}{6}\log\left(\frac{4\pi}{\nu_0(\s)} \right)
\end{equation}
We now claim that this condition can be satisfied by taking \(\varepsilon\) sufficiently small. Indeed, as noted in Remark \ref{rmq:Cepsilongdecroit}, the constant \(C > 0\) is non-increasing in \(\varepsilon\) and thus remains bounded away from zero as \(\varepsilon \to 0\), and consequently so does the right-hand side of \eqref{eq:diminuerepsilonpourcontract}, since $\nu_0(\s)=v\in(0,4\pi)$ is fixed. On the other hand, $\Lambda \to 0$ as $\varepsilon \to 0$. Indeed, since we are in dimension $2$, the curvature potential satisfies  
\[
\Delta_{\g} \xi_0 = e^{2\varphi}(R_0 - r),
\]  
where $e^{2\varphi}$ is the conformal factor of $g_0$ with respect to $\g$. A classical Schauder estimate combined with \eqref{eq:nearlyspherical} then implies that, up to a universal constant,  
\[
|\nabla^2 \xi_0| \lesssim \varepsilon.
\]  
Similarly, we clearly have, again up to a universal constant,  
\[
|(R_0 - r)g_0| \lesssim \varepsilon,
\]  
from which we deduce that \(|M_0| \to 0\) as \(\varepsilon \to 0\). Therefore, $\Lambda \underset{\varepsilon\to 0}{\longrightarrow} 0.$  

Letting $\varepsilon \to 0$ in \eqref{eq:diminuerepsilonpourcontract}, we conclude that there exists $\varepsilon_0 > 0$ such that \eqref{eq:diminuerepsilonpourcontract} holds for all $0 < \varepsilon < \varepsilon_0$, which completes the proof.\newline
 
\noindent
\textbf{Aknowledgements.} 
I would like to thank Marc Arnaudon and Michel Bonnefont for interesting discussions on the continuity equation of the Ricci flow, which provided valuable insights for section \ref{sec:KimMilman}. I would also like to thank Emanuel Milman for pointing out an inaccuracy in the formulation of Theorem \ref{thm:main} in an earlier version of this paper.

\bibliographystyle{plain}
\bibliography{mabibliographie}

\end{document}